\documentclass{colt2014} 

\usepackage{etoolbox}
\usepackage{algorithm2e}
\usepackage{graphicx}
\usepackage{url}
\usepackage{epstopdf}
\usepackage{color}

\newcommand{\RR}{\mathbb{R}}

\newcommand{\rank}{\mathrm{rank}}
\newcommand{\tr}{\mathrm{Tr}}

\title{Open problem: Tightness of maximum likelihood semidefinite relaxations}
\usepackage{times}

 \coltauthor{\Name{Afonso S. Bandeira} \Email{ajsb@math.princeon.edu}\\
 \addr Program in Applied and Computational Mathematics (PACM), Princeton University, Princeton NJ, USA
 \AND
 \Name{Yuehaw Khoo} \Email{ykhoo@princeton.edu}\\
 \addr Department of Physics, Princeton University, Princeton NJ, USA
 \AND
 \Name{Amit Singer} \Email{amits@math.princeton.edu}\\
 \addr Department of Mathematics and PACM, Princeton University, Princeton NJ, USA
}

\begin{document}

\maketitle

\begin{abstract}
We have observed an interesting, yet unexplained, phenomenon: Semidefinite programming (SDP) based relaxations of maximum likelihood estimators (MLE) tend to be tight in recovery problems with noisy data, even when MLE cannot exactly recover the ground truth. Several results establish tightness of SDP based relaxations in the regime where exact recovery from MLE is possible. However, to the best of our knowledge, their tightness is not understood beyond this regime. 
As an illustrative example, we focus on the generalized Procrustes problem.
\end{abstract}


\section{Introduction}

Recovery problems in many fields are commonly solved under the paradigm of maximum likelihood estimation. Despite the rich theory it enjoys, in many instances the parameter space is exponentially large and non-convex, often rendering the computation of the maximum likelihood estimator (MLE) intractable. It is then common to settle for heuristics, such as expectation-maximization. Unfortunately, popular iterative heuristics 
often get trapped in local minima and one usually does not know if the global optimum was achieved

A common alternative to these heuristics is the use of convex relaxations - to attempt optimizing (usually minimizing the minus log-likelihood) in a larger convex set that contains the parameter space of interest, as that allows one to leverage the power of convex optimization (e.g.~\cite{CandesTao_ConvexMatCompletion10}). The downside is that the solution obtained might not be in the original feasible set, forcing one to take an extra, potentially suboptimal, rounding step. The upside is that if it does lie in the original set then no rounding step is needed and one is guaranteed to have found the optimal solution to the original problem. Fortunately, this seems to often be the case in several problems.

We focus on a particular type of convex relaxations, namely semidefinite programming (SDP) based relaxations, and illustrate both these techniques and the open problem we wish to pose via the generalized orthogonal Procrustes problem (\cite{Nemirovski_sumsProcrustes}): There is an unknown underlying point cloud of $m$ points in $\RR^d$, $A\in\RR^{d\times m}$ (where columns of $A$ represent the coordinates of the points) and $n$ unknown orthogonal transformations in $\RR^d$, $\{O_i\}_{i = 1}^n$ satisfying $O_iO_i^T=I_{d\times d}$. We are given $n$ noisy measurements of the form 
\(A_i = O_iA + \sigma \eta\),
and, for simplicity, assume the noise to be white Gaussian, i.e., $\sigma>0$ and $\eta\in \RR^{d\times n}$ is a matrix with i.i.d. $\mathcal{N}(0,1)$ entries. The MLE for the unknowns $A$ and $\{O_i\}_{i = 1}^n$ is the minimizer of:
\(
\sum_{i=1}^n \left\|O_i^T A_i - A \right\|_F^2,
\)
constrained to $O_iO_i^T = I_{d\times d}$.
With the mild assumption that $\|A\|_F$ is fixed, the MLE is equivalent to an optimization problem only on the orthogonal transformations, we will refer to it as the quasi-MLE\footnote{If $\|A\|_F\neq 0$ is fixed, the MLE is the maximizer of $ \sum_i\tr\left(O_i^TA_i A^T\right)$, whose optimum $A$ is a multiple of $\sum_iO_i^TA_i$.}:
\begin{equation}\label{eq:Proscrustes:OrthogonalProblem}
\min_{O_1 , \dots, O_n \in \mathbf{O}_d}  \sum_{i<j} \left\|O_i^T A_i - O_j^T A_j \right\|_F^2.
\end{equation}
Unfortunately, the non-convexity and the exponential size of the search space renders problem (\ref{eq:Proscrustes:OrthogonalProblem}) intractable in general. 
Note that minimizing $\sum_{i<j} \left\|O_i^T A_i - O_j^T A_j \right\|_F^2$ has the same solution as maximizing $\sum_{i,j=1}^n \tr\left(A_jA_i^T O_iO_j^T \right)$ which is in turn, by taking $X = [O_1^T,\dots, O_n^T]^T[O_1^T,\dots, O_n^T]$ and $C\in\RR^{dn\times dn}$ with $d\times d$ blocks $C_{ij}^T = A_jA_i^T$, equivalent to
\begin{equation}\label{sdp:withrank}
\max \tr(CX) \quad \text{subject to } X\succeq 0,\ \forall_i\,X_{ii}=I_{d\times d},\ \rank(X)\leq d.
\end{equation}
The semidefinite relaxation in~\cite{Bandeira_LittleGrothendieckOd} can then be obtained by dropping the non-convex rank constraint, giving the following SDP
\begin{equation}\label{sdp:withoutrank}
\max \tr(CX) \quad \text{subject to } X\succeq 0,\ \forall_i\,X_{ii}=I_{d\times d}.
\end{equation}
When $d=1$ 
\cite{Abbe_Z2SynchER} showed that, below a certain level of outlier based noise, the solution of~(\ref{sdp:withoutrank}) achieves exact recovery (with high probability) and thus is a feasible point of~(\ref{sdp:withrank}). This is shown by constructing a dual certificate for the optimal point. 
However, for $d>1$, the $\ell_2$ nature of~(\ref{eq:Proscrustes:OrthogonalProblem}) together with the fact that there are infinitely many orthogonal transformations will render exact recovery under noise impossible, even by solving~(\ref{sdp:withrank}). Remarkably we observe that, even then, the relaxation~(\ref{sdp:withoutrank}) often manages to recover the quasi-MLE, the solution to~(\ref{sdp:withrank}).

This type of behavior has also been observed in the multireference alignment problem by~\cite{Bandeira_MultireferenceAlignment}, 
 in the global registration problem 
 by~\cite{Chaudhury_etal_GlobalRegistration}, and in camera motion estimation by~\cite{Ozyesil_CameraLocations_13}. Yet, to the best of our knowledge, there is no theoretical understanding of this {\it rank recovery} phenomenon. We note that there has been work on understanding the rank of solutions of random SDPs by~\cite{Amelunxen_rankrandomSDP} but the results hold only under specific distributions and do not apply to these problems. The difficulty of analyzing {\it rank recovery} lies in the fact that, unlike in exact recovery, we cannot identify the exact form of the MLE, rendering dual certificate arguments very difficult to carry out.

\section{An open problem}

\begin{figure}[h!]
  \centering
  \includegraphics[width=0.225\textwidth]{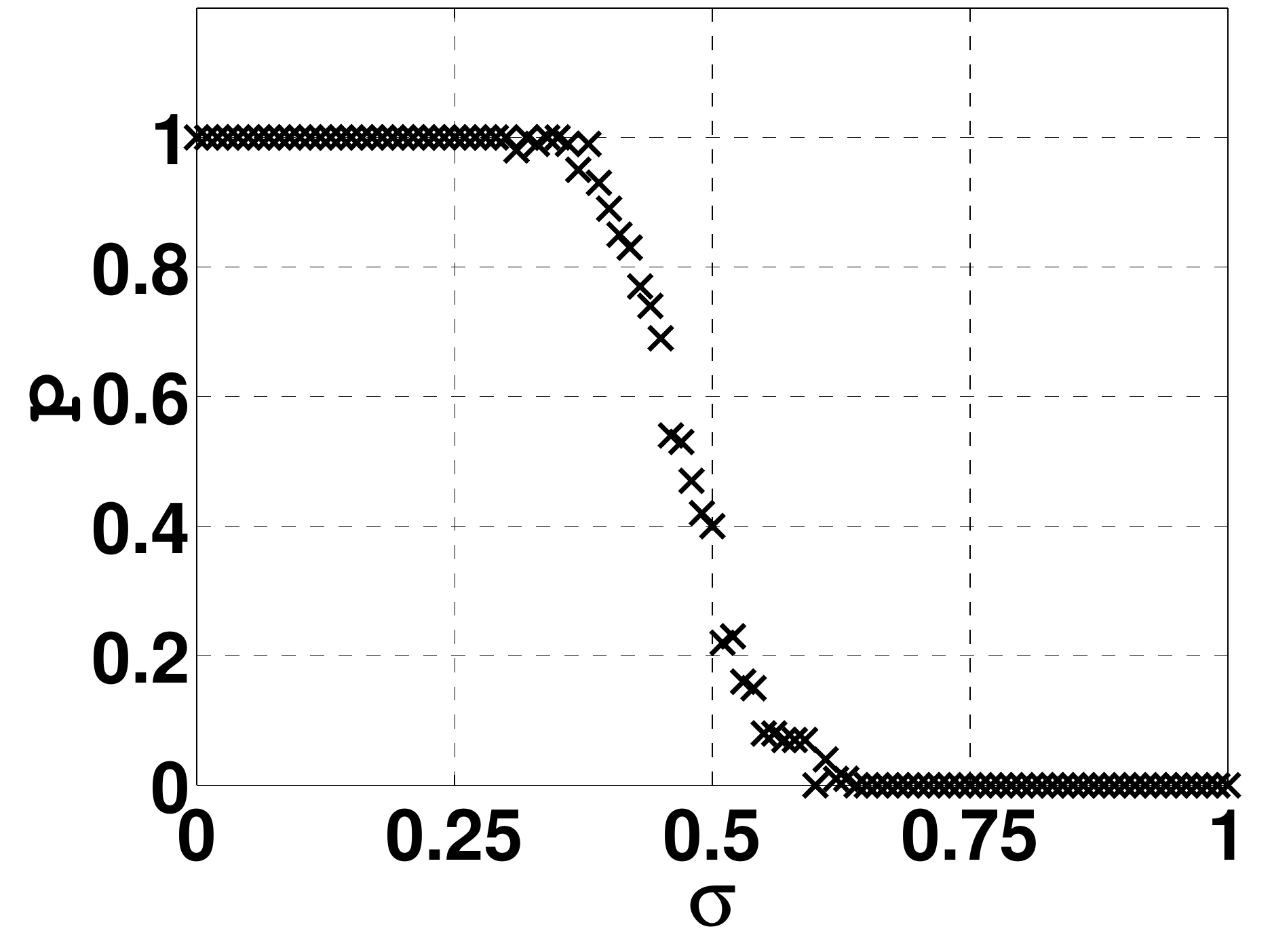} \quad \quad\includegraphics[width=0.225\textwidth]{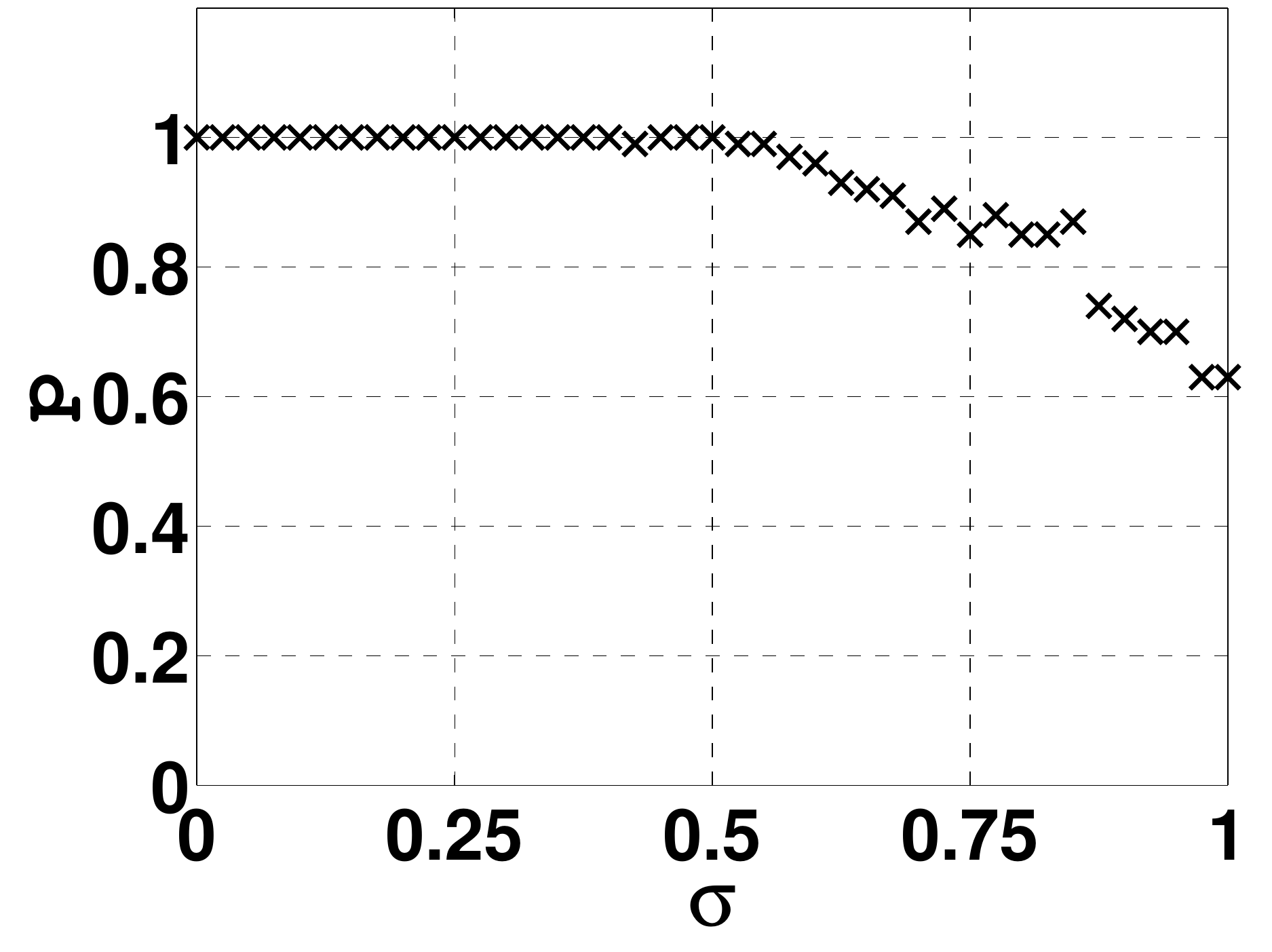}
    \caption{\footnotesize Fraction of trials (among 100) on which {\it rank recovery} was observed, for various values of noise $\sigma$. The left plot corresponds to the Procrustes problem (setting of Conjecture~\ref{Conjecture:Procrustes}) when $d=3, n = 20, m = 30$ and the right plot corresponds to the multi reference alignment problem treated in \cite{Bandeira_MultireferenceAlignment}. Both plots suggest that {\it rank recovery} happens with high probability, below certain noise levels.}\normalsize
    \label{fig:numericalevidence}
\end{figure}

Although we observe {\it rank recovery} in a variety of problems we formulate conjectures in two particularly simple problems, the generalized Procrustes problem~(Conjecture~\ref{Conjecture:Procrustes}) described above, and the multireference alignment problem treated in~\cite{Bandeira_MultireferenceAlignment}.
Numerical evidence supporting these conjecture is given in Figure~\ref{fig:numericalevidence}. We first pose the conjecture for generalized Procrustes.

\begin{conjecture}\label{Conjecture:Procrustes}
Let $d>1$, $n>2$, $m\geq d+1$, $A\in\RR^{d\times m}$ represent $m$ random points in $\RR^d$ (i.i.d. uniform random coordinates in $[0,1]^d$), and $\{O_i\}_{i=1}^n$ be any sequence of orthogonal transformations. Let $A_i = O_iA + \sigma \eta_i \in \mathbb{R}^{d\times m}$, where $\eta_i$ is a matrix with i.i.d. standard gaussian entries. There exists $\sigma_{\ast}>0$ such that, for $\sigma < \sigma_\ast$, with high probability, the solution of (\ref{sdp:withoutrank}) has rank $d$, hence matching the quasi-MLE (solution of (\ref{sdp:withrank})).
\end{conjecture}

The multireference alignment consists in estimating a $d$-dimensional signal by observing $n$ shifted noisy copies of it. For the sake of brevity, we will not describe the problem or the SDP based relaxation here (and refer the reader to~\cite{Bandeira_MultireferenceAlignment}) but take the opportunity to conjecture that a similar phenomenon happens: Below a certain noise level the solution to SDP-based relaxation in~\cite{Bandeira_MultireferenceAlignment} has rank $d$, thus matching the quasi-MLE. Although not going into details, we note that the SDP for this problem is considerably different than (\ref{sdp:withoutrank}), in particular, it has $\Omega\left(n^2d^2\right)$ positivity constraints. Also, this problem is discrete and so exact recovery is possible. However, exact recovery can be shown to be only possible for asymptotically vanishing levels of noise and here we conjecture {\it rank recovery} happens for a constant level of noise.

%


\bibliography{RankRecovery.bib}

\begin{thebibliography}{8}
\providecommand{\natexlab}[1]{#1}
\providecommand{\url}[1]{\texttt{#1}}
\expandafter\ifx\csname urlstyle\endcsname\relax
  \providecommand{\doi}[1]{doi: #1}\else
  \providecommand{\doi}{doi: \begingroup \urlstyle{rm}\Url}\fi

\bibitem[Abbe et~al.(2014)Abbe, Bandeira, Bracher, and Singer]{Abbe_Z2SynchER}
E.~Abbe, A.~S. Bandeira, A.~Bracher, and A.~Singer.
\newblock Linear inverse problems on {E}rd{\H{o}}s-{R}\'enyi graphs:
  Information-theoretic limits and efficient recovery.
\newblock \emph{IEEE International Symposium on Information Theory (ISIT2014)},
  to appear, 2014.

\bibitem[Amelunxen and B\"{u}rgisser(2014)]{Amelunxen_rankrandomSDP}
D.~Amelunxen and P.~B\"{u}rgisser.
\newblock Intrinsic volumes of symmetric cones and applications in convex
  programming.
\newblock \emph{Mathematical Programming}, pages 1--26, 2014.

\bibitem[Bandeira et~al.(2013)Bandeira, Kennedy, and
  Singer]{Bandeira_LittleGrothendieckOd}
A.~S. Bandeira, C.~Kennedy, and A.~Singer.
\newblock Approximating the little grothendieck problem over the orthogonal and
  unitary groups.
\newblock \emph{Available online at arXiv:1308.5207 [cs.DS]}, 2013.

\bibitem[Bandeira et~al.(2014)Bandeira, Charikar, Singer, and
  Zhu]{Bandeira_MultireferenceAlignment}
A.~S. Bandeira, M.~Charikar, A.~Singer, and A.~Zhu.
\newblock Multireference alignment using semidefinite programming.
\newblock \emph{5th Innovations in Theoretical Computer Science (ITCS 2014)},
  2014.

\bibitem[Candes and Tao(2010)]{CandesTao_ConvexMatCompletion10}
E.~J. Candes and T.~Tao.
\newblock The power of convex relaxation: Near-optimal matrix completion.
\newblock \emph{Information Theory, IEEE Transactions on}, 56\penalty0
  (5):\penalty0 2053--2080, May 2010.

\bibitem[Chaudhury et~al.(2013)Chaudhury, Khoo, and
  Singer]{Chaudhury_etal_GlobalRegistration}
K.~N. Chaudhury, Y.~Khoo, and A.~Singer.
\newblock Global registration of multiple point clouds using semidefinite
  programming.
\newblock \emph{arXiv:1306.5226 [cs.CV]}, 2013.

\bibitem[Nemirovski(2007)]{Nemirovski_sumsProcrustes}
A.~Nemirovski.
\newblock Sums of random symmetric matrices and quadratic optimization under
  orthogonality constraints.
\newblock \emph{Math. Program.}, 109\penalty0 (2-3):\penalty0 283--317, 2007.

\bibitem[Ozyesil et~al.(2013)Ozyesil, Singer, and
  Basri]{Ozyesil_CameraLocations_13}
O.~Ozyesil, A.~Singer, and R.~Basri.
\newblock Camera motion estimation by convex programming.
\newblock \emph{Available online at http://arxiv.org/abs/1312.5047 [cs.CV]},
  2013.

\end{thebibliography}

\end{document}